# Integers for Radical Extensions of Odd Prime Degree as Product of Subrings

Julius Kraemer

**Abstract.** For a radical extension $K$ of odd prime degree the ring $\mathcal{O}_K$ of integers is constructed as a product of subrings with the following property: for all prime divisors $q$ of the discriminant of $\mathcal{O}_K$ there is a $q$-maximal factor. The discriminant of $\mathcal{O}_K$ is the greatest common divisor of the discriminants of all factors. The results are applied to give a criterion for the monogeneity of $K$ where the opposite is not true.

The author is grateful to Keith Conrad and Hanson Smith for hints and references. The author is grateful to Guido Kraemer for proofreading a previous version.

## 1. Introduction

For $n \geq 2$ and $a \in \mathbb{Z}$, $a \neq \pm 1$, consider the polynomial $A(X) = X^n - a$ and assume that it is irreducible over $\mathbb{Z}$. Define $\alpha = \sqrt[n]{a}$, $\beta = \alpha - a$, $K = \mathbb{Q}(\alpha)$ and denote by $\mathcal{O}_K$ the ring of integers of $K$. Such radical extensions are also known as *pure extensions* in the literature.

The present paper is based on the findings of [3] where it has been characterized when $\mathcal{O}_K = \mathbb{Z}[\alpha]$. This is true under certain conditions (see Theorem 5.3 of [3]). Turning this theorem into the negative it follows that $\mathbb{Z}[\alpha] \subsetneq \mathcal{O}_K$ is equivalent to

  (i)   $a$ is not squarefree

or

  (ii)   There is a prime factor of $n$ which is a Wieferich prime to base $a$.

**Definition 1.1.** Let $q$ be a natural prime and $r \neq \pm 1$ a nonzero rational integer.
(i)   $q$ is a *Wieferich prime to base $r$* if $q^2$ divides $r^{q-1} - 1$.

(ii) If the context is clear we will speak of the *(non-)Wieferich case* if $q$ is (not) a Wieferich prime to base $r$.
(iii) Let $\mathcal{O} \subseteq \mathcal{O}_K$ be a subring. Then $\mathcal{O}$ is said to be *$q$-maximal* if the index $[\mathcal{O}_K : \mathcal{O}]$ is not divided by $q$.

For more information concerning Wieferich primes see Section 4 of [3] and the references cited there. Note that $q$ always divides $r^{q-1} - 1$ if $q$ and $r$ are coprime (Little Fermat Theorem). The definition of $q$-maximality can be found in 6.1.1 of [2].

In the present paper we assume that $n = p$ is an odd prime. The oddness of $p$ is not a real restriction because for $p = 2$ the ring of integers for radical extensions is already well-known.

In the main theorem we characterize $\mathcal{O}_K$ as a product of subrings of the form $\mathbb{Z}[\gamma_i]$, where $\mathbb{Z}[\gamma_i]$ is $q_i$-maximal and $q_i$ runs through all prime divisors of $a$; furthermore the minimal polynomial of all $\gamma_i$ is of the form $X^p - c_i$. In the Wieferich case an additional factor of the form $\mathbb{Z}[\beta']$ is necessary to ensure $p$-maximality. The discriminant of $\mathcal{O}_K$ is calculated as the greatest common divisor of the factors without using $\mathbb{Z}$ bases of $\mathcal{O}_K$.

The proof of the main theorem and its corollaries in Section 5 needs some preparation which is done in Sections 2 to 4:
In Section 2 preliminaries are handled which we need for the following sections.
In Section 3 the Wieferich case is handled: We prove that a specific subring of $\mathcal{O}_K$ is $p$-maximal.
In Section 4 the non-squarefree case is handled: For every prime factor $q$ of $a$ we construct a subring which is $q$-maximal.
In Section 6 we give examples and a criterion for the monogeneity of $K$ where the opposite is not true.

We denote the *greatest common divisor* of elements $x_1, \dots, x_n$ of a unique factorization domain by $(x_1, \dots, x_n)$ or $(x_i; 1 \leq i \leq n)$.

## 2. PRELIMINARIES

In this section we prove lemmas and propositions which are needed in the subsequent sections.

**Lemma 2.1.** *Let $q$ be a prime and $r, s$ rational integers. Then:*
(i) *The polynomial $X^q - r$ is irreducible over $\mathbb{Z}$ if and only if $X^q - rs^q$ is irreducible over $\mathbb{Z}$.*
(ii) *If $X^q - r$ is irreducible over $\mathbb{Z}$ then $\mathbb{Q}(\sqrt[q]{r}) = \mathbb{Q}(\sqrt[q]{rs^q})$.*
*Proof.* (i) A well-known theorem due to N. H. Abel (see Satz 277 together with Satz 180 of [7]) says that $X^q - r$ is irreducible over $\mathbb{Z}$ if and only if $r$ is not a $q$-th power in $\mathbb{Z}$. Then (i) follows immetiately.

(ii) This is immediate. □

**Remark 2.2.** As an immediate consequence of Lemma 2.1 we can assume in the following without restriction that in the prime decomposition $a = \prod_{i=1}^{m} q_i^{e_i}$ it holds that $1 \leq e_i \leq p-1$ for all $1 \leq i \leq m$. □

The following lemma introduces some further notation and properties of bases of Wieferich primes.

**Lemma 2.3.** *Let $r \neq \pm 1$ and $s$ be nonzero integers coprime to $p$. Then:*
(i) *If $p$ is not a Wieferich prime to base $r$ then, for $e \geq 1$, it holds that $p$ is a Wieferich prime to base $r^e$ if and only if $p$ divides $e$.*
(ii) *If $p$ is a Wieferich prime to base $r$ then $p$ is a Wieferich prime to base $r^e$ for all $e \geq 1$.*
(iii) *$p$ is a Wieferich prime to base $rs^p$ if and only if $p$ is a Wieferich prime to base $r$.*
   *If $s^p$ divides $r$ then $p$ is a Wieferich prime to base $\frac{r}{s^p}$ if and only if $p$ is a Wieferich prime to base $r$.*
(iv) *For $1 \leq j \leq m$ there are $1 \leq u_j \leq p-1$ and $v_j \geq 0$ such that*
   $e_j u_j - p v_j = 1$.
(v) *For $1 \leq j \leq m$ denote $c_j = \frac{a^{u_j}}{q_j^{p v_j}}$. Then $c_j$ is an integer where $q_j$ is a squarefree factor.*
   *$p$ is a Wieferich prime to base $c_j$ if and only if $p$ is a Wieferich prime to base $a$.*
(vi) $(c_j; 1 \leq j \leq m) = \prod_{j=1}^{m} q_j$.

*Proof.* (i) and (ii)  Decompose $r^{e(p-1)} - 1 = (r^{p-1} - 1) \cdot \sum_{i=0}^{e-1} r^{(p-1)i}$. Then (i) is clear because $p$ divides $r^{p-1} - 1$ (Little Fermat), and $p$ divides $\sum_{i=0}^{e-1} r^{(p-1)i}$ if and only if $p$ divides $e$ (again by Little Fermat). Also (ii) is clear because $p^2$ divides already $r^{p-1} - 1$.
(iii) Calculate
$$(rs^p)^{p-1} - 1 = r^{p-1} s^{p(p-1)} - 1 = s^{p(p-1)}(r^{p-1} - 1) + (s^{p(p-1)} - 1).$$
From (i) and (ii) it follows that $p^2$ divides $s^{p(p-1)} - 1$ hence, again applying (i) and (ii), it follows that $p^2$ divides $(rs^p)^{p-1} - 1$ if and only if $p^2$ divides $r^{p-1} - 1$ which is the first statement of (iii). The second statement of (iii) is proved analogously by using $\left(\frac{r}{s^p}\right)^{p-1} - 1 = \frac{1}{s^{p(p-1)}} \cdot \left((r^{p-1} - 1) - (s^{p(p-1)} - 1)\right)$ which is divided by $p^2$ if and only if the numerator is divided by $p^2$.
(iv) This follows immediately because $e_j$ and $p$ are coprime: Choose $u_j$ as a positive lift of $e_j^{-1} \bmod p$ to $\mathbb{Z}$ which lies between $1$ and $p-1$.
(v)  From (iv) it follows that
$$c_j = \frac{a^{u_j}}{q_j^{p v_j}} = \frac{\prod_{i=1}^{m} q_i^{e_i u_j}}{q_j^{p v_j}} = \frac{q_j^{e_j u_j} \cdot \prod_{i \neq j} q_i^{e_i u_j}}{q_j^{p v_j}} = q_j \cdot \prod_{i \neq j} q_i^{e_i u_j}$$

which is an integer where $q_j$ is a squarefree factor. The second statement follows from (i), (ii) and (iii) because $u_j$ is coprime to $p$ hence $p$ is a Wieferich prime to base $c_j$ if and only if $p$ is a Wieferich prime to base $a^{u_j}$ if and only if $p$ is a Wieferich prime to base $a$.
(vi)  This follows immediately from (v). □

**Remark 2.4.** *For $1 \leq j \leq m$ it holds that $e_j = 1$ if and only if $u_j = 1$ and $v_j = 0$. For $e_i \geq 2$ it follows that $u_j \geq 2$ and $v_j \geq 1$.*
*Proof.* This is immediate from Lemma 2.3 (iv). □

**Lemma 2.5.**
(i)  *It holds that $\mathbb{Z}[\alpha] = \mathbb{Z}[\beta]$.*
(ii) *The minimal polynomial of $\beta$ is*
$$B(X) = (X+a)^p - a = X^p + ap \cdot \left(\sum_{i=1}^{p-1} \frac{1}{p}\binom{p}{i} a^{p-1-i} X^i + \frac{a^{p-1}-1}{p}\right).$$
(iii) *Also $\{\beta^i; 0 \leq i \leq p-1\}$ is a power base of $\mathbb{Z}[\beta]$.*
(iv) *For a subring $\mathcal{O} \subseteq \mathcal{O}_K$ the discriminants of all $\mathbb{Z}$ bases of $\mathcal{O}$ are equal.*
(v)  *The discriminant $\mathrm{disc}(\mathbb{Z}[\beta])$ equals $(-1)^{\frac{p-1}{2}} \cdot a^{p-1} \cdot p^p$.*
(vi) *Exactly $p$ and the prime factors of $a$ can ramify in $\mathcal{O}_K$, in particular only $p$ and the prime factors of $a$ can divide $\mathrm{disc}(\mathcal{O}_K)$.*
*Proof.* (i)  This follows because $\beta \in \mathbb{Z}[\alpha]$ and $\alpha \in \mathbb{Z}[\beta]$.
(ii)  This holds because $A(X) = X^p - a$ is the minimal polynomial of $\alpha$.
(iii) This is clear.
(iv) This holds because the square of the determinant of every base change matrix $B$ between the power bases is $1$. Then the statement follows from the well-known theorem that the discriminants of two $\mathbb{Z}$ bases differ by the factor $\det(B)^2 = 1$ (see Proposition 1 in §2.7 of [8]).
(v)  For $0 \leq i, j \leq p - 1$ evaluate the trace of $\alpha^{i+j}$:
It is 0 for $i + j \neq p$ or $\neq 0$; it is $p$ if $i = j = 0$; it is $ap$ if $i + j = p$.
Then the discriminant can be calculated directly because it is the determinant of the trace matrix. Note that a diagonal matrix can be reached with $\frac{p-1}{2}$ interchanges of rows.
(vi) This is a well-known theorem, see Theorem 1 in §5.3 of [8]). □

We now state a criterion from which Dedekind's Index Criterion is a special case. It is stated in 6.1.4. of [2] and will be used intensively in the present paper. Denote by $\mathbb{F}_q$ the field with $q$ elements.

**Theorem 2.6** (Dedekind). *Let $K = \mathbb{Q}(\vartheta)$ be a number field, $T \in \mathbb{Z}[X]$ the minimal polynomial of $\vartheta$ and let q be a prime number.*
*Denote by $\overline{\phantom{x}}$ the reduction mod q (in $\mathbb{Z}$, $\mathbb{Z}[X]$ or $\mathbb{Z}[\vartheta]$). Let*
$$\overline{T}(X) = \prod_{i=1}^{k} \overline{T_i}(X)^{e_i}$$
*be the factorization of T mod q in $\mathbb{F}_q[X]$, and set*

$$G(X) = \prod_{i=1}^{k} T_i(X)$$

where the $T_i \in \mathbb{Z}[X]$ are arbitrary monic lifts of $\overline{T_i}$. Then:

(i) Let $H \in \mathbb{Z}[X]$ be a monic lift of $\overline{T(X)}/\overline{G(X)}$ and set

$$F = \frac{1}{q}(G \cdot H - T) \in \mathbb{Z}[X].$$

Then $\mathbb{Z}[\vartheta]$ is $q$-maximal if and only if $(\overline{F}, \overline{G}, \overline{H}) = 1$ in $\mathbb{F}_q[X]$.

(ii) Let $U$ be a monic lift of $\overline{T}/(\overline{F}, \overline{G}, \overline{H})$ to $\mathbb{Z}[X]$. Then $\mathcal{O}' = \mathbb{Z}[\vartheta] + \frac{1}{q}U(\vartheta)\mathbb{Z}[\vartheta]$ is a subring of $\mathcal{O}_K$, and if $m = \deg(\overline{F}, \overline{G}, \overline{H})$ then $[\mathcal{O}':\mathbb{Z}[\vartheta]] = q^m$ and $\mathrm{disc}(\mathcal{O}') = \frac{\mathrm{disc}(\mathbb{Z}[\vartheta])}{q^{2m}}$.

*Proof.* 6.1.4. of [2]. □

We apply Theorem 2.6 (i).

**Lemma 2.7.**
(i) Let $q$ be a prime factor of $a$. Then $\mathbb{Z}[\beta]$ is $q$-maximal if and only if $q$ does not occur squared in $a$.
(ii) The ring $\mathbb{Z}[\beta]$ is $p$-maximal if and only if $p$ is not a Wieferich prime to base $a$.

*Proof.* (i) Apply Theorem 2.6 (i) to $T = A$. Then $\overline{T} = X^p \bmod q$ hence $X$ and $X^{p-1}$ are monic lifts of $\overline{G}$ and $\overline{H}$, respectively. Then

$$F = \frac{1}{q} \cdot (G \cdot H - T) = \frac{a}{q} \not\equiv 0 \bmod q$$

if and only if $q$ is a squarefree factor of $a$. Because $(\overline{G}, \overline{H}) = (X, X^{p-1}) = X$ it follows that $(\overline{F}, \overline{G}, \overline{H}) = X$ if $q$ is not a squarefree factor of $a$, and $(\overline{F}, \overline{G}, \overline{H}) = 1$ if $q$ is a squarefree factor of $a$. This proves (i).

(ii) Apply Theorem 2.6 (i) to $T = B$. Then $\overline{T}(X) = X^p \bmod p$ hence $X$ and $X^{p-1}$ are again monic lifts of $\overline{G}$ and $\overline{H}$, respectively. Then

$$F = \frac{1}{p} \cdot (G \cdot H - T) = -a\left(\sum_{i=1}^{p-1} \frac{1}{p}\binom{p}{i} a^{p-1-i} X^i - \frac{a^{p-1}-1}{p}\right).$$

The constant factor of $F$ equals $0 \bmod p$ if and only if $p$ is a Wieferich prime to base $a$. Then it follows immediately that $(\overline{F}, \overline{G}, \overline{H}) = X$ if $p$ is a Wieferich prime to base $a$, and $(\overline{F}, \overline{G}, \overline{H}) = 1$ if $p$ is not a Wieferich prime to base $a$. This proves (ii). □

Next we apply Theorem 2.6 (ii).

**Lemma 2.8.** *Assume that $p$ is a Wieferich prime to base $a$ and let $U$ be a monic lift of $\overline{T}/(\overline{F}, \overline{G}, \overline{H})$ to $\mathbb{Z}[X]$. Then:*

INTEGERS FOR RADICAL EXTENSIONS OF ODD PRIME DEGREE AS PRODUCT OF SUBRINGS(i) It is possible to set $U(X) = X^{p-1}$ with $U(\beta) = \beta^{p-1}$. Then $\mathcal{O}' = \mathbb{Z}[\beta] + \frac{\beta^{p-1}}{p}\mathbb{Z}[\beta]$ is a ring and $\{\beta^i; 0 \leq i \leq p-2\} \cup \left\{\frac{\beta^{p-1}}{p}\right\}$ is a $\mathbb{Z} - base$ of $\mathcal{O}'$.

(ii) The degree of $(\overline{F}, \overline{G}, \overline{H})$ is 1, the index $[\mathcal{O}': \mathbb{Z}[\beta]]$ is $p^2$, and
$\text{disc}(\mathcal{O}') = (-1)^{\frac{p-1}{2}} \cdot p^{p-2} \cdot a^{p-1}$.

(iii) Denote $\beta' = \frac{\beta^{p-1}}{p}$. Then $\mathbb{Z}[\beta']$ is a subring of $\mathcal{O}_K$.

*Proof*. This follows immediately from Theorem 2.6 (ii). Note that $(\overline{F}, \overline{G}, \overline{H}) = X$ and apply Lemma 2.5 (v). Statement (iii) is also clear from (i). □

The following lemma has a more general context.

**Lemma 2.9.** *Let $L$ be a number field of degree $n$ and denote by $\mathcal{O}_L$ the ring of integers. Then*:

(i) Let $\mathcal{O}_1, \ldots, \mathcal{O}_\ell$ be subrings of $\mathcal{O}_L$. Then $\mathcal{O}' = \prod_{i=1}^{\ell} \mathcal{O}_i$ is a subring of $\mathcal{O}_L$ with $\mathcal{O}_i \subseteq \mathcal{O}'$ for all $1 \leq i \leq \ell$.

(ii) Let $\eta, \vartheta \in \mathcal{O}_L$ and assume that $\eta^2 \in \mathbb{Z}[\vartheta] + \eta\mathbb{Z}[\vartheta]$. Then
$$\mathbb{Z}[\eta] \cdot \mathbb{Z}[\vartheta] = \mathbb{Z}[\vartheta] + \eta\mathbb{Z}[\vartheta].$$

(iii) For a prime $q$ apply Theorem 2.6 to $\mathbb{Z}[\vartheta]$ and assume, with the notations of Theorem 2.6, that $(\overline{F}, \overline{G}, \overline{H}) \neq 1 \mod q$.
Denote $\eta = \frac{U(\vartheta)}{q}$. Then
$$\mathbb{Z}[\eta] \cdot \mathbb{Z}[\vartheta] = \mathbb{Z}[\vartheta] + \eta\mathbb{Z}[\vartheta].$$

(iv) In (i) let $\ell = 2$ and assume that the ranks of $\mathcal{O}_1$ and $\mathcal{O}_2$ are $n$. Then
$$\text{disc}(\mathcal{O}_2) = \text{disc}(\mathcal{O}_1) \cdot [\mathcal{O}_1 : \mathcal{O}_2]^2.$$
If $\mathcal{O}_2 \subseteq \mathcal{O}_1$ then $[\mathcal{O}_1 : \mathcal{O}_2]$ is an integer.

(v) In (iv) let $B$ be a base change matrix from $\mathcal{O}_1$ to $\mathcal{O}_2$. then
$$\det(B)^2 = [\mathcal{O}_1 : \mathcal{O}_2]^2.$$

(vi) In (i) it holds that $\text{disc}(\mathcal{O}')$ divides $(\text{disc}(\mathcal{O}_i); 1 \leq i \leq \ell)$.

*Proof*. The abelian group $\mathcal{O}'$ consists of all finite sums of products $\sum_{k=1}^{x} \sigma_{k1}\sigma_{k2}\ldots\sigma_{k\ell}$ with $\sigma_{ki} \in \mathcal{O}_i$.

(i) It is clear that $\mathcal{O}'$ is a ring with $\mathcal{O}_i \subseteq \mathcal{O}'$ for all $1 \leq i \leq \ell$ because $\sigma_i = \sigma_i \cdot 1$ with $\sigma_i \in \mathcal{O}_i$ and $1 \in \prod_{j \neq i} \mathcal{O}_j$.

(ii) "⊇": This is clear from (i) because both summands of the righthand side are contained in the lefthand side.
"⊆": For $0 \leq i, j \leq n-1$ it holds that $\eta^i \vartheta^j \in \mathbb{Z}[\eta] \cdot \mathbb{Z}[\vartheta]$. By induction on $i$ we show that $\eta^i \vartheta^j \in \mathbb{Z}[\vartheta] + \eta\mathbb{Z}[\vartheta]$ from which (ii) follows. For $i = 1$ the statement is clear. Assume that $\eta^i \vartheta^j \in \mathbb{Z}[\vartheta] + \eta\mathbb{Z}[\vartheta]$. Then there are $\varrho, \sigma \in \mathbb{Z}[\vartheta]$ such that $\eta^i \vartheta^j = \varrho + \eta\sigma$ hence $\eta^{i+1}\vartheta^j = \eta \cdot \eta^i \vartheta^j = \eta(\varrho + \eta\sigma) = \eta\varrho + \eta^2\sigma \in \mathbb{Z}[\vartheta] + \eta\mathbb{Z}[\vartheta]$ by assumption.

(iii) From Theorem 2.6 it follows that $\mathbb{Z}[\vartheta] + \eta\mathbb{Z}[\vartheta]$ is a ring which contains $\eta$. Then $\mathbb{Z}[\eta]$ is a subring of $\mathbb{Z}[\vartheta] + \eta\mathbb{Z}[\vartheta]$ hence (iii) follows from (ii).

(iv) This follows from (2.4) in III.2 of [4] because the prerequisites of page 94 of [4] are fulfilled: $\mathbb{Z}$ is a Dedekind domain, $\mathcal{O}_L$ has rank $n$ and this holds by assumption also for the subrings $\mathcal{O}_1$ and $\mathcal{O}_2$.

If $\mathcal{O}_2 \subseteq \mathcal{O}_1$ then $\mathcal{O}_1 \cdot \mathcal{O}_2 = \mathcal{O}_1$, and the statement follows from the Elementary Divisor Theorem (see Theorem 1 in §1.5 of [8]).
(v) Apply (iv) and the well-known statement $\operatorname{disc}(\mathcal{O}_2) = \operatorname{disc}(\mathcal{O}_1) \cdot \det(B)^2$ (see Proposition 1 in §2.7 of [8]).
(vi) This follows inductively from (iv) because $\mathcal{O}_i \subseteq \mathcal{O}'$ hence $[\mathcal{O}':\mathcal{O}_i]$ is an integer.
□

## 3. THE WIEFERICH CASE

In this section we assume that $p$ is a Wieferich prime to base $a$. Then $p$ does not divide $a$ and $\beta'$ is an integer in $K$ (Lemma 2.8 (iii)). We will prove

**Proposition 3.1.**
(i) If $p \geq 5$ then $\mathbb{Z}[\beta']$ is $p$-maximal.
(ii) If $p = 3$ then $\mathbb{Z}[\beta'] \cdot \mathbb{Z}[\alpha]$ is $p$-maximal.

The proof is done at the end of this section because several preparative statements are necessary.

**Remark 3.2** (matrix notation).
We will use matrix notation because we need the characteristic polynomial and the analysis modulo powers of $p$ seems to be easier.
Related to the base $\{\beta^i; 0 \leq i \leq p-1\}$ multiplication with $\beta$ is represented by the $p \times p$ matrix

$$M_\beta = \begin{pmatrix} 0 & 0 & \cdots & 0 & a - a^p \\ 1 & 0 & \cdots & 0 & -\binom{p}{1}a^{p-1} \\ 0 & 1 & \ddots & 0 & -\binom{p}{2}a^{p-2} \\ \vdots & \vdots & \ddots & \ddots & \vdots \\ 0 & 0 & \cdots & 1 & -\binom{p}{p-1}a \end{pmatrix}.$$

For details see Proposition 1 in §2.6 of [8]. The matrix $M_\beta$ has the following entries:
Row 1: $\quad m_{1p} = a - a^p; \quad m_{1j} = 0$ elsewhere
Row $i$ $(2 \leq i \leq p)$: $m_{i,i-1} = 1; \; m_{ip} = -\binom{p}{i-1}a^{p-(i-1)}; \; m_{ip} = 0$ elsewhere.
Multiplication with $\beta^k$ $(1 \leq k \leq p-1)$ is represented by $M_\beta^k = M_{\beta^k}$ hence multiplication with $\beta'$ is represented by $M_{\beta'} = \frac{1}{p}M_\beta^{p-1}$.
In $M_\beta$ every entry in the $p$-th column is divided by $p$ hence

$$M_\beta \equiv \begin{pmatrix} 0 & 0 & \cdots & 0 & 0 \\ 1 & 0 & \cdots & 0 & 0 \\ 0 & 1 & \ddots & 0 & 0 \\ \vdots & \vdots & \ddots & \ddots & \vdots \\ 0 & 0 & \cdots & 1 & 0 \end{pmatrix} \bmod p.$$

Then, for $1 \leq k \leq p - 1$,

$$M_\beta^k \equiv \begin{pmatrix} 0 & \cdots & \cdots & \cdots & \cdots & \cdots & 0 \\ \vdots & \vdots & \vdots & \vdots & \vdots & \vdots & \vdots \\ 0 & \cdots & \cdots & 0 & \cdots & \cdots & 0 \\ 1 & 0 & \cdots & \cdots & \cdots & \cdots & 0 \\ 0 & 1 & 0 & \cdots & \cdots & \cdots & 0 \\ \vdots & \ddots & \ddots & \ddots & \cdots & \cdots & 0 \\ 0 & \cdots & 0 & 1 & 0 & \cdots & 0 \end{pmatrix} \bmod p$$

with $1$ in the entries $(k + 1, 1), \dots, (p, p - k)$ (which is the $k$-th secondary diagonal below the main diagonal) and $0 \bmod p$ elsewhere. This means also that every entry in the columns $p - k + 1, \dots, p$ is divided by $p$.

For $1 \leq k \leq p - 1$ and $1 \leq i, j \leq p$ denote $M_\beta^k = \left(m_{ij}^{(k)}\right)$ and $m_{ij} = m_{ij}^{(1)}$.

This finishes Remark 3.2. ☐

**Lemma 3.3.** *Let $N = (n_{ij})$ be a $p \times p$ matrix with entries in the complex numbers and assume that $M_\beta \cdot N = N \cdot M_\beta$. Then $N' = (n'_{ij}) = M_\beta \cdot N$ has the following entries*:

(i) *The columns $2 \leq j \leq p$ of $N$ are shifted one column to the left*:
$n'_{ij} = n_{i,j+1}$ *for* $1 \leq i \leq p$, $1 \leq j \leq p - 1$.

(ii) *In the $p$ - th column of $N'$ the entry in the $i$ - th row is*
$n'_{1p} = (a - a^p) \cdot n_{pp}$
$n'_{ip} = n_{i-1,p} + \binom{p}{i-1} \cdot a^{p-(i-1)} \cdot n_{pp}$ *for* $2 \leq i \leq p$.

(iii) *If $p$ is a Wieferich prime to base $a$ then $p^2$ divides every entry in the first row of $M_\beta^k$ ($1 \leq k \leq p - 1$).*

*Proof.* (i) and (ii)   Use the structure of $M_\beta$ and matrix multiplication:
$N' = N \cdot M_\beta$ for (i), and $N' = M_\beta \cdot N$ for (ii).

(iii)   This follows from (i) and (ii) by induction on $k$, setting $N = M_\beta^k$, hence every entry in the first row is created by multiplication with $a - a^p$ or $0$ which both are divided by $p^2$. ☐

It is obvious that in the preceding lemma $M_\beta \cdot N = N \cdot M_\beta$ if $N$ is a power of $M_\beta$.

**Proposition 3.4.** *The minimal polynomial of $\beta'$ equals the characteristic polynomial of $\beta'$.*

*Proof.* Denote the $p \times p$ unit matrix by $E_p$. Then the characteristic polynomial $\chi$ of $\beta'$ is given by

$$\chi(X) = \det(X \cdot E_p - M_{\beta'}) = \det\left(X \cdot E_p - \frac{1}{p} M_\beta^{p-1}\right).$$

The matrix $M_\beta^{p-1}$ has the following properties (see Remark 3.2):
(a) It has integral entries in the columns $2, \ldots, p$ which are all divided by $p$.
(b) Every entry in the first row is divided by $p^2$ (Lemma 3.3 (iii)).
(c) The entry $(p, 1)$ equals $\frac{1}{p}$, all other entries in the first column are $0$.

Using the Laplace development for determinants along the first column it follows that $\chi$ is monic of degree $p$ and has integral coefficients. It remains to be shown that $\chi$ is the minimal polynomial of $\beta'$.

Because $\beta' \notin \mathbb{Q}$ it follows that the minimal polynomial $B'$ of $\beta'$ has degree $p$ as there are no further fields between $\mathbb{Q}$ and $K$. Our statement follows now from the well-known fact that $B'$ divides $\chi$ (which is the Cayley-Hamilton Theorem, see Satz 6 in Algebraische Ergänzung §2 of [1]). □

Next we analyze $M_\beta^k \mod p^2$.

**Proposition 3.5.** *For $1 \leq k \leq p - 1$ the following holds for $M_\beta^k$:*

(i) *For $1 \leq j \leq p - k$ the values are: $m_{ij}^{(k)} = 1$ for $j = i - k$ and $m_{ij}^{(k)} = 0$ for $j \neq i - k$.*

(ii) *In the first row the values are:*
$m_{1j}^{(k)} = 0$ *for* $1 \leq j \leq p - k$;
$m_{1,p+1-k}^{(k)} = a - a^p$;
$p(a - a^p)$ *divides* $m_{1j}^{(k)}$ *for* $p + 2 - k \leq j \leq p$,
*in particular* $p^3$ *divides* $p(a - a^p)$.

(iii) *For $2 \leq i \leq p$ it holds that $m_{ij}^{(k)} \equiv 0 \mod p^2$ for $p + i - k \leq j \leq p$.*

(iv) *For $2 \leq i \leq p$ it holds that*
$$m_{ij}^{(k)} \equiv -\binom{p}{i - (j - (p - k))} \cdot a^{p - (i - (j - (p - k)))} \mod p^2$$
*for $p - k + 1 \leq j \leq p - k + i - 1$.*

(v) *For $p \geq 5$ and $3 \leq k \leq p - 1$ it holds that $m_{2p}^{(k)} \not\equiv 0 \mod p^3$.*

*Proof.* (i) Consider in Remark 3.2 the first $p - k$ columns of $M_\beta^k$. The first $k$ rows are zeroes, the remaining $p - k$ rows form the unit matrix $E_{p-k}$. Then (i) follows immediately.

(ii) By Lemma 3.3 (i), applied inductively, the $p$-th column of $M_\beta$ is the $(p + 1 - k)$-th column of $M_\beta^k$ hence the first two statements of (ii) follow. For the third statement set $N = M_\beta$ in Lemma 3.3. Then $n'_{1p} = (a - a^p) \cdot n_{pp} = -(a - a^p)pa$, and the third statement follows inductively from Lemma 3.3 (i). It is immediate that $p^3$ divides $p(a - a^p)$ because we are in the Wieferich case.

(iii) Use again Lemma 3.3 (i) with $N = M_\beta$. From $p + i - k \leq j \leq p$ and $i \geq 2$ it follows that $k \geq i \geq 2$. The proof is done by induction on $k$. Put $k = 2$. Then $i = 2$ hence only row $2$ is relevant which means that $j = p$. Applying Lemma 3.3 (ii) yields that $n'_{2p} = n_{1p} + \binom{p}{1} a^{p-1} \cdot \binom{p}{p-1} a$ which is divided by $p^2$ because $n_{1p} = a - a^p$.
Now assume that $m_{ij}^{(k)} \equiv 0 \bmod p^2$ for $p + i - k \leq j \leq p$ and put $N = M_\beta^k$. Then $n'_{i+1,p} = n_{ip} + \binom{p}{i} a^{p-i} \cdot n_{pp}$. Because $n_{ip}$ is divided by $p^2$ (induction hypothesis) and $n_{pp}$ is divided by $p$, the statement follows for the $p$-th column. Then statement (iii) follows by applying again Lemma 3.3 (i).

<u>(iv)</u> The proof is done by induction on $k$. For $k = 1$ it follows that $p \leq j \leq p - 2 + i$ hence $j = p$. Then
$$-\binom{p}{i - (j - (p-k))} \cdot a^{p - (i - (j - (p-k)))} \equiv -\binom{p}{i-1} \cdot a^{p-(i-1)} \equiv m_{ip} \bmod p^2$$
for $2 \leq i \leq p$ hence (iv) holds for $k = 1$. Now assume that (iv) holds for $k$. Put again $N = M_\beta^k$ and $N' = M_\beta \cdot N = M_\beta^{k+1}$. Then

$$\begin{aligned} m_{ip}^{(k+1)} &= n'_{ip} = n_{i-1,p} + \binom{p}{i-1} \cdot a^{p-(i-1)} \cdot n_{pp} \equiv n_{i-1,p} \equiv m_{i-1,p}^{(k)} \\ &\equiv -\binom{p}{(i-1) - (p - (p-k))} \cdot a^{p - ((i-1) - (p - (p-k)))} \\ &\equiv -\binom{p}{i - (k+1)} \cdot a^{p-(i-(k+1))} \bmod p^2 \text{ for } i \geq k+1. \end{aligned}$$

(For the first congruence use Lemma 3.3 (iii) and $n_{pp} \equiv 0 \bmod p$; for the third congruence use the induction hypothesis and the results for $k$.) Hence statement (iv) follows for the $p$-th column. Then statement (iv) follows inductively from Lemma 3.3 (i).

<u>(v)</u> From (iv) it follows, with $j = p$, that $m_{ip}^{(k)} \equiv -\binom{p}{i-k} \cdot a^{p-(i-k)} \bmod p^2$ for $1 \leq k \leq i - 1$, in particular $m_{ip}^{(k)} \equiv 0 \bmod p$ and $\not\equiv 0 \bmod p^2$. From (iii) it follows that $m_{2p}^{(k)} \equiv 0 \bmod p^2$ (put $i = 2$ and $j = p$, then $p + 2 - k \leq p$ hence $k \geq 2$).
For $k = 2$ it follows, with $N = M_\beta$, from Lemma 3.2 (iii) that
$$\begin{aligned} m_{1p}^{(2)} &= -(a - a^p)\binom{p}{p-1} a \\ m_{2p}^{(2)} &= (a - a^p) + p^2 a^p \equiv 0 \bmod p^2 \\ m_{ip}^{(2)} &= m_{i-1,p} - \binom{p}{i-1} \cdot a^{p-(i-1)} \cdot \binom{p}{p-1} \cdot a \\ &= -\binom{p}{i-2} \cdot a^{p-(i-2)} - p\binom{p}{i-1} \cdot a^{p+2-i} \end{aligned}$$
which is $\equiv 0 \bmod p$ and $\not\equiv 0 \bmod p^2$ for $3 \leq i \leq p$.
For $k = 3$ (this is possible because $p \geq 5$) it follows, with $N = M_\beta^2$, again from Lemma 3.3 (ii), applied twice, that
$$m_{2p}^{(3)} = m_{1p}^{(2)} + pa^{p-1} m_{pp}^{(2)} = -(a - a^p) p \cdot a - p^3 \cdot a^{p+1} - p\binom{p}{p-2} \cdot a^{p+1}$$
which is $\equiv 0 \bmod p^2$ and $\not\equiv 0 \bmod p^3$;
$$m_{ip}^{(3)} \equiv -\binom{p}{i-3} \cdot a^{p-(i-3)} \bmod p^2 \text{ for } i \geq 4$$
which is $\equiv 0 \bmod p$ and $\not\equiv 0 \bmod p^2$ for $3 \leq i \leq p$.
The proof is now done by induction on $3 \leq k \leq p - 2$. Put $N = M_\beta^k$ and assume that $m_{2p}^{(k)} \not\equiv 0 \bmod p^3$. Then $m_{2p}^{(k+1)} = m_{1p}^{(k)} + pa^{p-1} m_{pp}^{(k)}$.

From $m_{pp}^{(k)} \equiv 0 \bmod p$ and $\not\equiv 0 \bmod p^2$ (by (iv), because $k \leq p-1 < p$) it follows that $p \cdot a^{p-1} \cdot m_{pp}^{(k)} \equiv 0 \bmod p^2$ and $\not\equiv 0 \bmod p^3$. Then statement (v) follows because $m_{1p}^{(k)}$ is (by (ii)) divided by $p(a - a^p)$ hence $m_{1p}^{(k)} \equiv 0 \bmod p^3$. □

The reason why $p = 3$ has to be treated separately is illustrated in Remark and Example 3.10.

**Corollary 3.6.**
(i) For $j = p + (i-1) - k$ and $i \leq k+1$ it holds that $m_{ij}^{(k)} \equiv -p \bmod p^2$.
(ii) For $k = p - 1$ Proposition 3.5 reads as follows:
   (i) First column: $m_{p1}^{(p-1)} = 1$; $m_{i1}^{(p-1)} = 0$ for $1 \leq i \leq p-1$.
   (ii) First row: $m_{11}^{(p-1)} = 0$; $m_{12}^{(p-1)} = a - a^p$;
      $p(a - a^p)$ divides $m_{1j}^{(p-1)}$ for $3 \leq j \leq p$,
      in particular $p(a - a^p) \equiv 0 \bmod p^3$.
   (iii) For $2 \leq i \leq p$ it holds that $m_{ij}^{(p-1)} \equiv 0 \bmod p^2$ for $i+1 \leq j \leq p$.
      In particular there are integers $m'_{ij}$ such that $m_{ij}^{(p-1)} = p^2 \cdot m'_{ij}$.
   (iv) For $2 \leq i \leq p$ it holds that
      $m_{ij}^{(p-1)} \equiv -\binom{p}{i-(j-1)} \cdot a^{p-(i-(j-1))} \bmod p^2$ for $2 \leq j \leq i$.
      In particular $m_{ij}^{(p-1)} = p \cdot m'_{ij}$ where $m'_{ij}$ is a unit mod $p^2$.
(iii) For $k = p - 1$ it holds that $m_{ii}^{(p-1)} \equiv -p \bmod p^2$ for $2 \leq i \leq p$.
(iv) The matrix $M_\beta^{p-1}$ reads $\bmod\, p^2$ as follows:

$$M_\beta^{p-1} \equiv \begin{pmatrix} 0 & 0 & 0 & \cdots & \cdots & 0 & 0 \\ 0 & -p & 0 & \cdots & \cdots & 0 & 0 \\ 0 & -\binom{p}{2}a^{p-2} & -p & 0 & \cdots & 0 & 0 \\ 0 & -\binom{p}{3}a^{p-3} & -\binom{p}{2}a^{p-2} & -p & \ddots & \vdots & \vdots \\ \vdots & \vdots & \vdots & \ddots & \ddots & \ddots & \vdots \\ 0 & -\binom{p}{p-2}a^2 & -\binom{p}{p-3}a^3 & \cdots & \ddots & -p & 0 \\ 1 & -\binom{p}{p-1}a & -\binom{p}{p-2}a^2 & \cdots & \cdots & -\binom{p}{2}a^{p-2} & -p \end{pmatrix}.$$

*Proof.* (i) From the assumption it follows that
$$i - (j - (p-k)) = i - (p + (i-1) - k - (p-k)) = 1$$
hence $m_{ij}^{(k)} \equiv -p \cdot a^{p-1} \bmod p^2$ from Proposition 3.5 (iv). Then (i) follows because we are in the Wieferich case.
(ii) to (iv) They follow immediately by applying Proposition 3.5 for $k = p - 1$. □

Now it is possible to calculate the minimal polynomial $\chi(X)$ of $\beta' \bmod p$ (see Proposition 3.4).

**Proposition 3.7.**
(i) $\chi(X) \equiv X(X+1)^{p-1} \bmod p$
(ii) $\mathbb{Z}[\beta']$ is $p$-maximal if and only if $p^2$ does not divide $\chi(-1)$.

*Proof.* (i) $\chi(X) = \det(X \cdot E_p - M_{\beta'}) = \det\left(X \cdot E_p - \frac{1}{p} \cdot M_\beta^{p-1}\right)$

$$= \det \begin{pmatrix} X & -\frac{1}{p}m_{12}^{(p-1)} & \cdots & -\frac{1}{p}m_{1p}^{(p-1)} \\ 0 & X - \frac{1}{p}m_{22}^{(p-1)} & \cdots & -\frac{1}{p}m_{2p}^{(p-1)} \\ \vdots & \vdots & \ddots & \vdots \\ -\frac{1}{p} & \cdots & \cdots & X - \frac{1}{p}m_{pp}^{(p-1)} \end{pmatrix}$$

$$= X \cdot \det \begin{pmatrix} X - \frac{1}{p}m_{22}^{(p-1)} & \cdots & -\frac{1}{p}m_{2p}^{(p-1)} \\ \vdots & \ddots & \vdots \\ -\frac{1}{p}m_{p2}^{(p-1)} & \cdots & X - \frac{1}{p}m_{pp}^{(p-1)} \end{pmatrix}$$

$$-\frac{1}{p} \cdot \det \begin{pmatrix} -\frac{1}{p}m_{12}^{(p-1)} & \cdots & \cdots & -\frac{1}{p}m_{1p}^{(p-1)} \\ X - \frac{1}{p}m_{22}^{(p-1)} & \ddots & \cdots & -\frac{1}{p}m_{2p}^{(p-1)} \\ \vdots & \ddots & \ddots & \vdots \\ -\frac{1}{p}m_{p-1,2}^{(p-1)} & \cdots & X - \frac{1}{p}m_{p-1,p-1}^{(p-1)} & -\frac{1}{p}m_{p-1,p}^{(p-1)} \end{pmatrix}.$$

Look at the first summand: All entries $-\frac{1}{p}m_{ij}^{(p-1)}$ above the main diagonal are congruent $0 \bmod p$ (Corollary 3.6 (iv)). Because $m_{ii}^{(p-1)} \equiv -p$ for $2 \leq i \leq p$ it follows that the first summand is $\bmod p$ congruent to $X(X+1)^{p-1}$.
Look at the second summand: It equals

$$\det\begin{pmatrix} \frac{1}{p^2}m_{12}^{(p-1)} & \cdots & & \cdots & \frac{1}{p^2}m_{1p}^{(p-1)} \\ X - \frac{1}{p}m_{22}^{(p-1)} & \ddots & & \ddots & -\frac{1}{p}m_{2p}^{(p-1)} \\ \vdots & \ddots & & \ddots & \vdots \\ -\frac{1}{p}m_{p-1,2}^{(p-1)} & \cdots & X - \frac{1}{p}m_{p-1,p-1}^{(p-1)} & & -\frac{1}{p}m_{p-1,p}^{(p-1)} \end{pmatrix}$$

because $-\frac{1}{p}$ can be put in the first row. The right column is congruent $0 \bmod p$ because $p^2$ divides every entry from the second row on (Corollary 3.6 (iv)) and $p(a - a^p)$ divides the entry $m_{1p}^{(p-1)}$ (Corollary 3.6 (ii) property (ii)) because we are in the Wieferich case. Then the second summand is congruent $0 \bmod p$ which proves (i).

<u>(ii)</u>   Apply Theorem 2.6 to $T = \chi$. Then $X \cdot (X + 1)$ and $(X + 1)^{p-2}$ are monic lifts of $\overline{G}$ and $\overline{H}$, respectively, and $F = \frac{1}{p}(X \cdot (X + 1)^{p-1} - \chi(X))$ with $F(-1) = \frac{-\chi(-1)}{p}$ and $(\overline{G}, \overline{H}) = X + \overline{1}$. It follows that $(\overline{F}, \overline{G}, \overline{H}) = X + \overline{1}$ if and only if $X + \overline{1}$ divides $\overline{F}$, and this is equivalent to $\frac{-\chi(-1)}{p} \equiv 0 \bmod p$. Then (ii) follows because $\frac{-\chi(-1)}{p} \not\equiv 0 \bmod p$ if and only if $-\chi(-1) \not\equiv 0 \bmod p^2$. □

Before showing that $p^2$ does not divide $\chi(-1)$ for $p \geq 5$ we need the following

**Lemma 3.8.** *Let $r \geq 2$ be an integer and $N = (n_{ij})$ be a $r \times r$ matrix with integer entries and the following property: There is an integer $s$ such that $n_{ij} \equiv 0 \bmod s$ for $i \leq j$. Then*

$$\det(N) \equiv n_{1r} \cdot \prod_{i=2}^{r} n_{i,i-1} \bmod s^2.$$

*Proof.* By assumption there are integers $n'_{ij}$ with $n_{ij} = s \cdot n'_{ij}$ for $i \leq j$. Then (standard calculation for determinants)

$$\det(N) = \sum_{i=1}^{r} (-1)^{i-1} \cdot s \cdot n'_{1i} \cdot \det(N_i)$$

where the $(r - 1) \times (r - 1)$ matrix $N_i$ is defined by deleting the first row and the $i$-th column from $N$. For $1 \leq i \leq r - 1$ the entries in the right column of $N_i$ are $s \cdot n'_{ij}$ hence $s$ divides $\det(N_i)$, and then $s^2$ divides every summand of $\det(N)$ except the $r$-th one. Then $\det(N) \equiv s \cdot n'_{1r} \cdot \det(N_r) \bmod s^2$. The entries on the main diagonal of $N_r$ are $n_{i,i-1}$ ($2 \leq i \leq r$), above the main diagonal all entries are divided by $s$ hence $\det(N_r) \equiv \prod_{i=2}^{r} n_{i,i-1} \bmod s$ because $N_r$ is $\bmod s$ a lower triangular matrix. This proves the statement. □

**Remark 3.9.**  $\det(N) \equiv 0 \bmod s^2$ *if and only if $s^2$ divides $n_{1r}$ or $s$ divides one of the $n_{i,i-1}$ ($2 \leq i \leq r$).* □

Now we can prove Proposition 3.1.

*Proof of Proposition* 3.1.

<u>(i)</u>   From the prerequisite we have $p \geq 5$.
From Corollary 3.6 (ii) (property (iii)), (iii) and (iv) it follows that all entries of $M_\beta^{p-1}$ above the main diagonal are congruent $0 \bmod p^2$, and all entries below the main diagonal are not divided by $p^2$.
Now look at the first summand of $\chi(X)$ in the proof of Proposition 3.7 (i). Because $m_{ii}^{(p-1)} \equiv -p \bmod p^2$ (Corollary 3.6 (iii)) the first summand reads $\bmod\, p^2$ as follows:

$$X \cdot \det\begin{pmatrix} X+1 & \ddots & pm'_{2p} \\ \vdots & \ddots & \vdots \\ m'_{p2} & \cdots & X+1 \end{pmatrix}$$

where the $m'_{ij}$ ($2 \leq j < i \leq p$) below the main diagonal are units $\bmod\, p$. Evaluation at $-1$ gives

$$-\det\begin{pmatrix} 0 & \ddots & pm'_{2p} \\ \vdots & \ddots & \vdots \\ m'_{p2} & \cdots & 0 \end{pmatrix}$$

which is by Lemma 3.8 (with $r = s = p$) congruent to $-pm'_{2p} \cdot \prod_{i=3}^{p} m'_{i,i-1} \bmod p^2$. This is not congruent $0 \bmod p^2$ because all factors except $p$ are units $mod\, p$, and then also $mod\, p^2$: $m'_{2p}$ by Proposition 3.5 (v), and $m'_{i,i-1}$ by Corollary 3.6 (ii) property (iv). Note that $m'_{2p}$ is a unit $\bmod\, p^2$ because $m_{2p}^{(p-1)} = p^2 m'_{2p}$ is not divided by $p^3$.

Now look at the second summand of $\chi(X)$ in the proof of Proposition 3.7 (ii) and evaluate at $-1$. Because we are in the Wieferich case it follows from Corollary 3.6 (ii) property (ii) that $p^3$ divides $m_{1j}^{(p-1)}$ for $3 \leq j \leq p$. Then the second summand reads $\bmod\, p^2$ as follows:

$$\det \begin{pmatrix} \dfrac{a-a^p}{p^2} & pm'_{13} & \cdots & \cdots & \cdots & \cdots & pm'_{1p} \\ 0 & pm'_{23} & \ddots & \cdots & \cdots & \cdots & pm'_{2p} \\ \vdots & \ddots & \ddots & \ddots & \cdots & \cdots & \vdots \\ m'_{i2} & \cdots & 0 & \ddots & \ddots & pm'_{ij} & \vdots \\ \vdots & \vdots & \vdots & \ddots & \ddots & \ddots & \vdots \\ m'_{p-2,2} & m'_{p-2,3} & \vdots & \vdots & \ddots & \ddots & pm'_{p-2,p} \\ m'_{p-1,2} & m'_{p-1,3} & \cdots & \cdots & \cdots & 0 & pm'_{p-1,p} \end{pmatrix}$$

which is congruent $0 \mod p^2$ because all entries in the last two columns are divided by $p$.

Putting together the results for both summands statement (i) follows for $p \geq 5$.

(ii) From the prerequisite we have $p = 3$.

From Lemma 2.9 (iii) it follows that $\mathbb{Z}[\beta'] \cdot \mathbb{Z}[\alpha] = \mathbb{Z}[\alpha] + \beta'\mathbb{Z}[\alpha]$. From $\operatorname{disc}(\mathbb{Z}[\alpha]) = -27a^2$ (Lemma 2.5 (v)) and Theorem 2.6 (ii) it follows that $\operatorname{disc}(\mathbb{Z}[\beta'] \cdot \mathbb{Z}[\alpha]) = -3a^2$ which is not divided by $9$. Then (ii) follows from Proposition 1 in §2.7 of [8].

The proof of Proposition 3.1 is now completed. $\square$

**Remark and Example 3.10.** For $p = 3$ we see from the proof of Proposition 3.5 (v) that $m^{(2)}_{2p} \equiv 0 \mod p^2$ but not necessarily $\not\equiv 0 \mod p^3$. For $p = 3$ the entry $m^{(2)}_{23} = (a - a^3) + 9a^3 = a(8a^2 + 1)$ in the matrix $M^2_\beta$ can be congruent $0 \mod 27$. Put, for example, $a = 19$. Then $m^{(2)}_{23} = 27 \cdot 19 \cdot 107$.

The minimal polynomial of $\beta'$ is
$$X^3 - a^2 X^2 + a^2 \frac{a^2 + 2}{3} X - a^2 \frac{(a^2 - 1)^2}{27}.$$
Applying Theorem 2.6 to $\mathbb{Z}[\beta']$ shows that $\mathbb{Z}[\beta']$ is 3-maximal if and only if $\frac{a^2-1}{9} \not\equiv 1 \mod 3$. Again $a = 19$ is an example where $\mathbb{Z}[\beta']$ is not 3-maximal. $\square$

## 4. The Non-Squarefree Case

Due to Remark 2.2 it holds for the exponents of the prime decomposition of $a = \prod_{i=1}^m q_i^{e_i}$ without restriction that $1 \leq e_i \leq p - 1$ for all $1 \leq i \leq m$.

Let $1 \leq i \leq m$. For $\alpha_i = \sqrt[p]{q_i}$ it follows that $\alpha = \prod_{i=1}^m \alpha_i^{e_i}$. For $c_i = \frac{a^{u_i}}{q_i^{pv_i}}$ (introduced in Lemma 2.3 (v)) denote $\gamma_i = \frac{\alpha^{u_i}}{q_i^{v_i}}$. Then $\gamma_i = \sqrt[p]{c_i} = \alpha_i \cdot \prod_{\substack{k=1,\\ k \neq i}}^m \alpha_k^{e_k u_i}$ (proof of Lemma 2.3 (v)).

**Proposition 4.1.** *Let* $1 \leq i \leq m$. *Then*:
(i) $\gamma_i \in \mathcal{O}_K$ *and* $\mathbb{Q}(\gamma_i)$ *is a radical extension.*
(ii) $\mathbb{Z}[\gamma_i]$ *is* $q_i$-*maximal.*
(iii) $\operatorname{disc}(\mathbb{Z}[\gamma_i]) = (-1)^{\frac{p-1}{2}} \cdot p^p \cdot c_i^{p-1}$.
(iv) $\mathbb{Z}[\gamma_i]$ *is* $p$-*maximal if and only if* $\mathbb{Z}[\alpha]$ *is* $p$-*maximal.*
(v) *The ring* $\mathcal{O}^* = \prod_{j=1}^{m} \mathbb{Z}[\gamma_i]$ *is* $q_i$-*maximal for all* $1 \leq i \leq m$.
(vi) $\mathbb{Z}[\alpha] \subseteq \mathcal{O}^*$.

*Proof.* (i) The number $\gamma_i = \frac{\alpha^{u_i}}{q_i^{v_i}}$ is a root of $X^p - c_i \in \mathbb{Z}[X]$ and lies in $K$. Then (i) follows.
(ii) This follows from the representation $c_i = q_i \cdot \prod_{\substack{k=1 \\ k \neq i}}^{m} q_k^{e_k u_i}$ (see Lemma 2.3, proof of (v)) where $q_i$ is a squarefree factor. Now (ii) follows from Lemma 2.7 (i).
(iii) The argumentation is identical to Lemma 2.5 (v) because $\{\gamma_i^k; 0 \leq k \leq p-1\}$ is a power base of $\mathbb{Z}[\gamma_i]$ and $X^p - c_i$ is the minimal polynomial of $\gamma_i$.
(iv) Apply Lemma 2.3 (v) and Lemma 2.7 (ii).
(v) This follows from Lemma 2.9 (i).
(vi) This follows from (v) and the prime decomposition of $a$. □

## 5. MAIN THEOREM AND COROLLARIES

We can start immediately with the main theorem of this paper.

**Theorem 5.1.** *With the notations from the preceding sections the following holds*:
(i) *If* $p$ *is a Wieferich prime to base* $a$ *then* $\mathcal{O}_K = \mathbb{Z}[\beta'] \cdot \mathcal{O}^*$.
(ii) *If* $p$ *is not a Wieferich prime to base* $a$ *then* $\mathcal{O}_K = \mathcal{O}^*$.
(iii) *The discriminant* $\operatorname{disc}(\mathcal{O}_K)$ *equals* $(-1)^{\frac{p-1}{2}} \cdot p^x \cdot \prod_{j=1}^{m} q_j^{p-1}$
 *with* $x = p - 2$ *in statement* (i) *and* $x = p$ *in statement* (ii).
(iv) $\operatorname{disc}(\mathcal{O}_K) = \begin{cases} \left(\operatorname{disc}(\mathbb{Z}[\gamma_j]); 1 \leq j \leq m\right) \text{ in the non-Wieferich case.} \\ \left(\operatorname{disc}(\mathbb{Z}[\beta']), \operatorname{disc}(\mathbb{Z}[\gamma_j]); 1 \leq j \leq m\right) \text{ in the Wieferich case.} \end{cases}$

*Proof.* By Lemma 2.9 (i) $\mathcal{O}^*$ is a subring of $\mathcal{O}_K$, every $\mathbb{Z}[\gamma_j]$ is contained in $\mathcal{O}^*$ and $\mathcal{O}^*$ is $q_j$-maximal for all $1 \leq j \leq m$ by Proposition 4.1 (v). As consequence from Lemma 2.5 (vi) it remains to show the $p$-maximality of $\mathcal{O}^*$ and $\mathbb{Z}[\beta'] \cdot \mathcal{O}^*$, respectively.
(i) If $\mathbb{Z}[\alpha]$ is not $p$-maximal then, by Proposition 3.1, $\mathbb{Z}[\beta']$ and $\mathbb{Z}[\beta'] \cdot \mathbb{Z}[\alpha]$ are $p$-maximal and 3-maximal, respectively. From $\mathbb{Z}[\alpha] \subseteq \mathcal{O}^*$ (Proposition 4.1 (vi)) it follows now that $\mathbb{Z}[\beta'] \cdot \mathcal{O}^*$ is $p$-maximal which proves (i).
(ii) Let $\mathbb{Z}[\alpha]$ be $p$-maximal. Then, by Proposition 4.1 (iv), $\mathbb{Z}[\gamma_j]$ is $p$-maximal for all $1 \leq j \leq m$ hence $\mathcal{O}^*$ is also $p$-maximal which proves (ii).
(iii) and (iv) From Proposition 4.1 (ii) and, in the Wieferich case, Theorem 2.6 (ii) it follows immediately that $\left(\operatorname{disc}(\mathbb{Z}[\gamma_j]); 1 \leq j \leq m\right) = (-1)^{\frac{p-1}{2}} \cdot p^p \cdot \prod_{j=1}^{m} q_j^{p-1}$ and $(\mathrm{disc}(\mathbb{Z}[\beta']), \mathrm{disc}(\mathbb{Z}[\gamma_j]); 1 \le j \le m) = (-1)^{\frac{p-1}{2}} \cdot p^{p-2} \cdot \prod_{j=1}^{m} q_j^{p-1}$. It remains to be shown that these numbers equal $\mathrm{disc}(\mathcal{O}_K)$.

Firstly assume the non-Wieferich case. Then $\mathcal{O}_K = \mathcal{O}^*$. From Lemma 2.9 (vi) we know that $\mathrm{disc}(\mathcal{O}^*)$ divides $(\mathrm{disc}(\mathbb{Z}[\gamma_i]); 1 \le i \le m)$. We also know that the only divisors of $\mathrm{disc}(\mathcal{O}^*)$ are $p, q_1, \dots, q_m$ (Lemma 2.5 (vi) and Proposition 4.1 (iii)).

Assume that there is $1 \le j \le m$ such that $q_j$ divides $\frac{(\mathrm{disc}(\mathbb{Z}[\gamma_i]); 1 \le i \le m)}{\mathrm{disc}(\mathcal{O}^*)}$. Then $q_j$ divides every $\frac{\mathrm{disc}(\mathbb{Z}[\gamma_i])}{\mathrm{disc}(\mathcal{O}^*)}$ $(1 \le i \le m)$, in particular $q_j$ divides $\frac{\mathrm{disc}(\mathbb{Z}[\gamma_j])}{\mathrm{disc}(\mathcal{O}^*)} = [\mathcal{O}^*:\mathbb{Z}[\gamma_j]]^2$ by Lemma 2.9 (iv) which is a contradiction to the $q_j$-maximality of $\mathbb{Z}[\gamma_j]$. With the same argument it is shown that $p$ does not divide $\frac{(\mathrm{disc}(\mathbb{Z}[\gamma_i]); 1 \le i \le m)}{\mathrm{disc}(\mathcal{O}^*)}$ because all $\mathbb{Z}[\gamma_i]$ are $p$-maximal by Proposition 4.1 (iv).

Now assume the Wieferich case. Because $\mathbb{Z}[\beta']$ is idempotent (which means $\mathbb{Z}[\beta']^2 = \mathbb{Z}[\beta']$) it follows that $\mathcal{O}_K = \mathbb{Z}[\beta'] \cdot \mathcal{O}^* = \prod_{j=1}^{m}(\mathbb{Z}[\beta'] \cdot \mathbb{Z}[\gamma_j])$ where every factor $\mathbb{Z}[\beta'] \cdot \mathbb{Z}[\gamma_j]$ is $p$-maximal and $q_j$-maximal. Furthermore $\mathbb{Z}[\beta'] \cdot \mathbb{Z}[\gamma_j] = \mathbb{Z}[\gamma_j] + \beta'\mathbb{Z}[\gamma_j]$ by Lemma 2.9 (iii) with $\mathrm{disc}(\mathbb{Z}[\beta'] \cdot \mathbb{Z}[\gamma_j]) = \frac{\mathrm{disc}(\mathbb{Z}[\gamma_j])}{p^2}$ by Theorem 2.6 (ii). The statement in the Wieferich case now follows with the same argument as in the non-Wieferich case replacing $\mathbb{Z}[\gamma_j]$ by $\mathbb{Z}[\beta'] \cdot \mathbb{Z}[\gamma_j]$ $(1 \le j \le m)$. The proof of (iii) and (iv) is now complete. □

Statement (iii) of Theorem 5.1 is well-known (Section 3 of [10] and a consequence of Theorem 1.1 of [5]). The proof given here is different from these approaches.

The first corollary shows that it is possible to integrate the factors of $\mathcal{O}_K$ if certain prerequisites are fulfilled.

### Corollary 5.2.
(i) For $1 \le i, j \le m$ with $i \ne j$ assume that $e_i = e_j$. Then
$$\mathbb{Z}[\gamma_i] \cdot \mathbb{Z}[\gamma_j] = \mathbb{Z}\left[\frac{\alpha^{u_i}}{q_i^{v_i} \cdot q_j^{v_j}}\right]$$
and $\mathbb{Z}\left[\frac{\alpha^{u_i}}{q_i^{v_i} \cdot q_j^{v_j}}\right]$ is $q_i$- and $q_j$-maximal.

(ii) For $1 \le i \le m$ denote $a_i = \prod_{\substack{j=1 \\ e_j=e_i}}^{m} q_j$ and $\mathcal{O}_i = \prod_{\substack{j=1 \\ e_j=e_i}}^{m} \mathbb{Z}[\gamma_j]$. Then
$$\mathcal{O}_i = \mathbb{Z}\left[\frac{\alpha^{u_i}}{\prod_{\substack{j=1 \\ e_j=e_i}}^{m} q_j^{v_j}}\right].$$

(iii) Let $\ell$ be the number of different $e_i$ for $1 \le i \le m$. Then
$\mathcal{O}^* = \prod_{k=1}^{\ell} \mathcal{O}_k$ where $k$ runs over all different $e_j$.

*Proof.* (i) From $e_i = e_j$ it follows that $u_i = u_j$ because $1 \le u_i \le p-1$. Then $\gamma_i = \frac{\alpha^{u_i}}{q_i^{v_i}} = \frac{\alpha^{u_i}}{q_i^{v_i} \cdot q_j^{v_j}} \cdot q_j^{v_j}$ lies in the righthand side. The same argument applies to $\gamma_j$

INTEGERS FOR RADICAL EXTENSIONS OF ODD PRIME DEGREE AS PRODUCT OF SUBRINGShence "⊆" follows. Because $q_i$ and $q_j$ are coprime there are integers $r$ and $s$ with
$1 = q_i^{v_i} r + q_j^{v_j} s$ hence $\frac{\alpha^{u_i}}{q_i^{v_i} \cdot q_j^{v_j}} = \frac{\alpha^{u_i} \cdot (q_i^{v_i} r + q_j^{v_j} s)}{q_i^{v_i} \cdot q_j^{v_j}} = \frac{\alpha^{u_i}}{q_j^{v_j}} r + \frac{\alpha^{u_i}}{q_i^{v_i}} s$ which lies in the lefthand side, and then also "⊇" follows. An analogous argument as in Proposition 4.1 (v) shows that $\mathbb{Z}\left[\frac{\alpha^{u_i}}{q_i^{v_i} \cdot q_j^{v_j}}\right]$ is $q_i$- and $q_j$-maximal.

<u>(ii)</u>  This is shown inductively from (i) because all distinct $q_i, q_j$ with $e_i = e_j$ are pairwise coprime.

<u>(iii)</u>  This follows because all exponents of all prime factors of $a$ are covered by exactly one $\mathcal{O}_k$.  □

Generally it is enough to assume in Corollary 5.2 that $e_i \equiv e_j \bmod p$ but due to Remark 2.2 it is possible to assume that $1 \leq e_i \leq p - 1$ hence $e_i = e_j$.

The next corollary shows how, in the Wieferich case, $\mathbb{Z}[\beta']$ can be replaced by other rings.

**Corollary 5.3.** *In the Wieferich case* $\mathbb{Z}[\beta']$ *can be replaced in Theorem* 5.1 *by any of the rings* $\mathbb{Z}\left[\frac{(\gamma_i - c_i)^{p-1}}{p}\right] (1 \leq i \leq m)$.

*Proof.* The polynomial $X^p - c_i$ is the minimal one for $\gamma_i$ hence, by Lemma 2.3 (v), the situation for $c_i$ is the same as for $a$. From Proposition 3.1 it follows that also all rings $\mathbb{Z}\left[\frac{(\gamma_i - c_i)^{p-1}}{p}\right]$ and $\mathbb{Z}\left[\frac{(\gamma_i - c_i)^2}{3}\right] \cdot \mathbb{Z}[\alpha]$ are $p$-maximal ($p \geq 5$) and 3-maximal subrings of $\mathcal{O}_K$, respectively.  □

## 6. EXAMPLES AND A CRITERION FOR MONOGENEITY

Firstly, we construct $\mathbb{Z}$ bases of $\mathcal{O}_K$.

**Proposition 6.1.** *With the notations from above the following holds*:
(i)  For $0 \leq k \leq p - 1$ and $1 \leq j \leq m$ there are $t_{kj} \geq 0$ and 
$0 \leq e'_{kj} \leq p - 1$ such that $ke_i = pt_{ki} + e'_{ki}$ and $\frac{\alpha^k}{\prod_{j=1}^m q_j^{t_{kj}}} = \prod_{j=1}^m \alpha_j^{e'_{kj}} \in \mathcal{O}_K$.
If $k = 0$ then $t_{kj} = e'_{kj} = 0$ else $e'_{kj} \geq 1$.
(ii) For $1 \leq j \leq m$ and $1 \leq k \leq p - 2$ it holds that $t_{k+1,j} \geq t_{kj}$, and for $1 \leq k \leq p - 1$ it holds that $(p - 1 - k)e_j \geq t_{p-1,j} - t_{kj}$.
(iii) If $p$ is not a Wieferich prime to base $a$ then
$$B' = \left\{\frac{\alpha^k}{\prod_{j=1}^m q_j^{t_{kj}}}; \ 0 \leq k \leq p - 1\right\}$$
*is a* $\mathbb{Z}$ *base of* $\mathcal{O}_K$.

(iv) *If $p$ is a Wieferich prime to base $a$ then*

$$B'' = \left\{\frac{\alpha^k}{\prod_{j=1}^m q_j^{t_{kj}}}; \ 0 \leq k \leq p-2\right\} \cup \left\{\frac{(\alpha-a)^{p-1}}{p \cdot \prod_{j=1}^m q_j^{t_{p-1,j}}}\right\}$$

*is a $\mathbb{Z}$ base of $\mathcal{O}_K$.*

We omit the (lengthy but straightforward) proof because $\mathbb{Z}$ bases of $\mathcal{O}_K$ have already been constructed in [10] and, as a special case, in [6]. The methods used in our proof are quite similar to the one used in [10].

The following example illustrates the results developed here.

**Example 6.2.**

(i) Put $p = 5$ and let $q_1, \ldots, q_5$ be different primes. Set $a' = q_1 q_2^2 q_3^4 q_4^7 q_5^5$ and assume that $5$ is not a Wieferich prime to base $a'$. From Remark 2.2 and Lemma 2.3 (iii) it follows that $\mathcal{O}_K$ is completely determined by $a = q_1 q_2^2 q_3^4 q_4^2$. Use the notations from above. Then:

$m = 4$;
$e_1 = 1, \ e_2 = 2, \ e_3 = 4, \ e_4 = 2$;
$\alpha = \sqrt[5]{a}, \ \alpha_i = \sqrt[5]{q_i}$ with $\alpha = \alpha_1 \alpha_2^2 \alpha_3^4 \alpha_4^2$;
$(u_1, v_1) = (1, 0), \ (u_2, v_2) = (3, 1), \ (u_3, v_3) = (4, 3), \ (u_4, v_4) = (3, 1)$;
$c_1 = a, \ c_2 = \frac{a^3}{q_2^5}, \ c_3 = \frac{a^4}{q_3^{15}}, \ c_4 = \frac{a^3}{q_4^5}$;
$\gamma_i = \sqrt[5]{c_i} \ (1 \leq i \leq 4)$ with $\gamma_1 = \alpha, \ \gamma_2 = \frac{\alpha^3}{q_2}, \ \gamma_3 = \frac{\alpha^4}{q_3^3}, \ \gamma_4 = \frac{\alpha^3}{q_4}$;
$\mathcal{O}_K = \mathbb{Z}[\alpha] \cdot \mathbb{Z}\left[\frac{\alpha^3}{q_2}\right] \cdot \mathbb{Z}\left[\frac{\alpha^4}{q_3^3}\right] \cdot \mathbb{Z}\left[\frac{\alpha^3}{q_4}\right]$ (Theorem 5.1 (i)) with
$\mathbb{Z}\left[\frac{\alpha^3}{q_2}\right] \cdot \mathbb{Z}\left[\frac{\alpha^3}{q_4}\right] = \mathbb{Z}\left[\frac{\alpha^3}{q_2 \cdot q_4}\right]$ (Corollary 5.2 (i)).

Next we compute a $\mathbb{Z}$ base of $\mathcal{O}_K$.

For $0 \leq k \leq 4$ the $ke_i = pt_{ki} + e'_{ki}$ from Proposition 6.1 (i) can be represented in the following matrix with entries $(t_{ki}, e'_{ki})$ (5 columns numerated from $0$ to $4$, and $4$ rows numerated from $1$ to $4$):

$$\begin{pmatrix}
 & k=0 & k=1 & k=2 & k=3 & k=4 \\
e_1 = 1 & (0,0) & (0,1) & (0,2) & (0,3) & (0,4) \\
e_2 = 2 & (0,0) & (0,2) & (0,4) & (1,1) & (1,3) \\
e_3 = 4 & (0,0) & (0,4) & (1,3) & (2,2) & (3,1) \\
e_4 = 2 & (0,0) & (0,2) & (0,4) & (1,1) & (1,3)
\end{pmatrix}.$$

Then $\prod_{i=1}^4 q_i^{t_{1i}} = 1$, $\prod_{i=1}^4 q_i^{t_{2i}} = q_3$, $\prod_{i=1}^4 q_i^{t_{3i}} = q_2 q_3^2 q_4$, $\prod_{i=1}^4 q_i^{t_{4i}} = q_2 q_3^3 q_4$ hence $B' = \left\{1, \ \alpha, \ \frac{\alpha^2}{q_3}, \ \frac{\alpha^3}{q_2 q_3^2 q_4}, \ \frac{\alpha^4}{q_2 q_3^3 q_4}\right\}$ is a $\mathbb{Z}$ base of $\mathcal{O}_K$.

(ii) Assume in (i) that $5$ is a Wieferich prime to base $a$. Then $\beta' = \frac{(\alpha-a)^4}{5}$ lies in $\mathcal{O}_K$ and $\mathbb{Z}[\beta']$ is 5-maximal (Proposition 3.1 (i)). For $\mathcal{O}^* = \mathbb{Z}[\alpha] \cdot \mathbb{Z}\left[\frac{\alpha^3}{q_2}\right] \cdot \mathbb{Z}\left[\frac{\alpha^4}{q_3^3}\right] \cdot \mathbb{Z}\left[\frac{\alpha^3}{q_4}\right]$ it follows that $\mathcal{O}_K = \mathbb{Z}[\beta'] \cdot \mathcal{O}^*$, and a $\mathbb{Z}$ base of $\mathcal{O}_K$ is $B'' = \left\{1, \alpha, \frac{\alpha^2}{q_3}, \frac{\alpha^3}{q_2 q_3^2 q_4}, \frac{(\alpha-a)^4}{5 q_2 q_3^3 q_4}\right\}$.

(iii) Put $p = 11$ and $a = 3^2 = 9$. Then $11$ is a Wieferich prime to base $3$ because $11^2$ divides $3^{10} - 1$. Due to Lemma 2.3 (ii) $11$ is also a Wieferich prime to base $9$. Then $\alpha = \sqrt[11]{9}$ and $\alpha_1 = \sqrt[11]{3}$, $m = 1$, $e_1 = 2$, $(u_1, v_1) = (6,1)$; $c_1 = \frac{9^6}{3^{11}} = 3$, $\gamma_1 = \sqrt[11]{3}$; $\beta' = \frac{(\sqrt[11]{9}-9)^{10}}{11}$ with $\mathcal{O}_K = \mathbb{Z}[\beta'] \cdot \mathbb{Z}[\sqrt[11]{3}]$.

For $0 \leq k \leq 10$ the number $ke_1 = 2k = pt_{k1} + e'_{k1}$ from Proposition 6.1 (i) can be represented as follows:

$$\begin{pmatrix} k & 0 & 1 & 2 & 3 & 4 & 5 \\ (t_{k1}, e'_{k1}) & (0,0) & (0,2) & (0,4) & (0,6) & (0,8) & (0,10) \end{pmatrix}$$

$$\begin{pmatrix} k & 6 & 7 & 8 & 9 & 10 \\ (t_{k1}, e'_{k1}) & (1,1) & (1,3) & (1,5) & (1,7) & (1,9) \end{pmatrix}.$$

Then a $\mathbb{Z}$ base of $\mathcal{O}_K$ is
$B'' = \left\{1, \sqrt[11]{9}, \sqrt[11]{9}^2, \ldots, \sqrt[11]{9}^5, \frac{\sqrt[11]{9}^6}{3}, \ldots, \frac{\sqrt[11]{9}^9}{3}, \frac{(\sqrt[11]{9}-9)^{10}}{11 \cdot 3}\right\} = \left\{\sqrt[11]{3}^i \; (0 \leq i \leq 9), \frac{(\sqrt[11]{9}-9)^{10}}{33}\right\}.$
  □

There is an application of Theorem 5.1 if all exponents of the prime decomposition of $a$ have the same congruence $\mod p$ and $\mathbb{Z}[\alpha]$ is $p$-maximal.

**Definition 6.3.** A number field $L$ of degree $n$ over the rational numbers is called *monogenic* if there is a $\vartheta$ in the ring of integers $\mathcal{O}_L$ such that $\{1, \vartheta, \ldots, \vartheta^{n-1}\}$ is a power base of $\mathcal{O}_L$.

The literature concerning (non-)monogeneity of number fields is extensive; for a starting point see Section 1 of [9] and the references cited there.

**Proposition 6.4.** *If $\mathbb{Z}[\alpha]$ is $p$-maximal and all exponents of the prime decomposition of $a$ which are not divided by $p$ have the same congruence $\mod p$ then $K$ is monogenic with $\mathcal{O}_K = \mathbb{Z}\left[\frac{\alpha^{u_1}}{\prod_{i=1}^m q_i^{v_i}}\right]$.*

*Proof.* This is an immediate consequence from Remark 2.2 and Corollary 5.2 (iii) because $\ell = 1$ and $\mathbb{Z}[\alpha]$ is $p$-maximal.  □

In particular this proposition can be applied if $a$ is squarefree which is Theorem 2.1 of [3] in the case of prime degree. In this case $u_i = 1$ and $v_j = 0$ in Lemma 2.3 (iv) hence $c_j = a$ which yields $\mathcal{O}_K = \mathbb{Z}[\alpha]$ if $\mathbb{Z}[\alpha]$ is $p$-maximal.

The other direction of this proposition is not true as the following example shows. It extends Exercise 10B to Chapter V of [8].

**Example 6.5.** Let $p = 3$ and $q_1, q_2$ be distinct primes. Put $a = q_1 q_2^2$. Then, with the notations of this paper, $\alpha_1 = \sqrt[3]{q_1}$, $\alpha_2 = \sqrt[3]{q_2}$ and $\alpha = \sqrt[3]{a} = \alpha_1 \alpha_2^2$. Additionally denote $\alpha' = \alpha_1^2 \alpha_2$. Then, from Theorem 5.1 (iii) and Proposition 6.1 (iii) and (iv), $\mathcal{O}_K$ has $\mathbb{Z}$ base $\{1, \alpha, \alpha'\}$ with $\mathrm{disc}(\mathcal{O}_K) = -27(q_1 q_2)^2$ in the non-Wieferich case and $\mathbb{Z}$ base $\{1, \alpha, \frac{(\alpha-a)^2}{3q_2}\}$ with $\mathrm{disc}(\mathcal{O}_K) = -3(q_1 q_2)^2$ in the Wieferich case (where $3$ does not divide $a$), respectively. Note that the prerequisites of Proposition 6.4 are not fulfilled.

Let $\eta = t_0 + t_1 \alpha + t_2 \alpha' \in \mathcal{O}_K$ with $t_0, t_1, t_2 \in \mathbb{Z}$. In the non-Wieferich case $\eta$ is an arbitrary element from $\mathcal{O}_K$, in the Wieferich case there are further elements in $\mathcal{O}_K$.

Firstly, we assume the non-Wieferich case. Using $\alpha^2 = q_2 \alpha'$, $\alpha'^2 = q_1 \alpha$ and $\alpha \alpha' = q_1 q_2$ it follows that $\eta^2 = t_0^2 + 2 t_1 t_2 q_1 q_2 + (t_2^2 q_1 + 2 t_0 t_1) \alpha + (t_1^2 q_2 + 2 t_0 t_2) \alpha'$ hence the base change matrix from $\{1, \alpha, \alpha'\}$ to $\{1, \eta, \eta^2\}$ is

$$\begin{pmatrix} 1 & 0 & 0 \\ t_0 & t_1 & t_2 \\ t_0^2 + 2 t_1 t_2 q_1 q_2 & t_2^2 q_1 + 2 t_0 t_1 & t_1^2 q_2 + 2 t_0 t_2 \end{pmatrix} \text{ with determinant}$$

$t_1^3 q_2 - t_2^3 q_1$. Then $\mathrm{disc}(\mathbb{Z}[\eta]) = -27(q_1 q_2)^2 \cdot (t_1^3 q_2 - t_2^3 q_1)^2$ (see Proposition 1 in §2.7 of [8]) which corrects statement (e) of the mentioned exercise. From $\mathrm{disc}(\mathcal{O}_K) = -27(q_1 q_2)^2$ it follows that that $\mathcal{O}_K = \mathbb{Z}[\eta]$ if and only if there is a solution $(t_1, t_2)$ such that $(t_1^3 q_2 - t_2^3 q_1)^2 = 1$. Solutions for $(q_1, q_2, t_1, t_2)$ are, for example, $x_1 = (2, 17, 1, 2)$, $x_2 = (5, 41, 1, 2)$, $x_3 = (11, 19, -5, -6)$. Because $3$ is a Wieferich prime to base $5 \cdot 41^2$ and not in the other cases it follows that $\mathcal{O}_K = \mathbb{Z}[\eta]$ for $(q_1, q_2) = (2, 17)$ and $(11, 19)$ but not for $(q_1, q_2) = (5, 41)$.

Now we assume the Wieferich case. Denote $\beta'' = \frac{(\alpha-a)^2}{3q_2}$ and note that $\frac{4+5a^2}{9}, \frac{1+8a^2}{9}$ are integers because we are in the Wieferich case. Also $\frac{1+2a^2}{3}$ is an integer because $3$ and $a$ are coprime. Then straightforward calculations show:

$\alpha^2 = -a^2 + 2a\alpha + 3q_2 \beta''$

$\alpha \beta'' = q_1 q_2 \frac{1+2a^2}{3} - q_1 q_2 a \alpha - 2a \beta''$

$\beta''^2 = -q_1 a \frac{4+5a^2}{9} + q_1 \frac{1+8a^2}{9} \alpha + 2 q_1 q_2 a \beta''$.

Now let $\vartheta = s_0 + s_1 \alpha + s_2 \beta''$ ($s_0, s_1, s_2 \in \mathbb{Z}$) be an arbitrary element from $\mathcal{O}_K$. Then $\vartheta^2 = s_0' + s_1' \alpha + s_2' \beta''$ with

$s_0' = s_0^2 - s_1^2 a^2 - s_2^2 q_1 a \frac{4+5a^2}{9} + 2 s_1 s_2 q_1 q_2 \frac{1+2a^2}{3}$,

$s_1' = 2 s_1^2 a + s_2^2 q_1 \frac{1+8a^2}{9} + 2 s_0 s_1 - 2 s_1 s_2 q_1 q_2 a$,

$s_2' = 3 s_1^2 q_2 + 2 s_2^2 q_1 q_2 a + 2 s_0 s_2 - 4 s_1 s_2 a$.

Then the base change matrix from $\{1, \alpha, \beta''\}$ to $\{1, \vartheta, \vartheta^2\}$ is

$$\begin{pmatrix} 1 & 0 & 0 \\ s_0 & s_1 & s_2 \\ s_0' & s_1' & s_2' \end{pmatrix} \text{ with determinant}$$

$$\Delta = s_1 s_2' - s_2 s_1' = 3s_1^3 q_2 + 4s_1 s_2^2 q_1 q_2 a - 6s_1^2 s_2 a - s_2^3 q_1 \frac{1+8a^2}{9}$$

hence $9\Delta = (3s_1 - 2s_2 q_1 q_2)^3 q_2 - s_2^3 q_1$. With the same argument as in the non-Wieferich case it follows that that $\mathcal{O}_K = \mathbb{Z}[\vartheta]$ if and only if there is a solution $(s_1, s_2)$ such that $\Delta = 1$. A solution is, for example, $(q_1, q_2, s_1, s_2) = (2, 11, 15, 1)$ as is easily calculated. □

The results of Example 6.5 can be summarized as follows:

**Remark 6.6.** *In the situation of Example* 6.5 *the following holds*:
(i)   *If* 3 *is not a Wieferich prime to base* $a$ *then* $K$ *is monogenic if and only if the equation* $t_1^3 q_2 - t_2^3 q_1 = 1$ *has a solution with integers* $t_1, t_2$. *If* $K$ *is monogenic then* $\mathcal{O}_K = \mathbb{Z}[t_1 \alpha + t_2 \alpha']$.
(ii)  *If* 3 *is a Wieferich prime to base* $a$ *then* $K$ *is monogenic if and only if the equation* $(3s_1 - 2s_2 q_1 q_2)^3 q_2 - s_2^3 q_1 = 9$ *has a solution with integers* $s_1, s_2$. *If* $K$ *is monogenic then* $\mathcal{O}_K = \mathbb{Z}\left[s_1 \alpha + s_2 \frac{(\alpha-a)^2}{3q_2}\right]$. □

RUWIDO-BOGEN 11
D-85635 HÖHENKIRCHEN-SIEGERTSBRUNN, GERMANY

Julius_Kraemer@gmx.de

Mathematisches Institut der Ludwig-Maximilians-Universität München from 1972 until 1988